\documentclass[12pt]{article}

\usepackage{amsmath}
\usepackage{amssymb}

\normalbaselines
\newtheorem{Theorem}{Theorem}[section]

\newtheorem{Proposition}[Theorem]{Proposition}
\newtheorem{Lemma}[Theorem]{Lemma}
\newtheorem{Corollary}[Theorem]{Corollary}
\newtheorem{Remark}[Theorem]{Remark}

\newcommand{\rl}{\mathbb{R}}

\newcommand{\dt}{\delta}
\newcommand{\ep}{\varepsilon}
\newcommand{\ld}{\lambda}
\newcommand{\mm}{{\cal M}}

\newcommand{\nb}{\nabla}

\newcommand{\df}{\displaystyle\frac}

\begin{document}

\def\qed{\hbox{\hskip 6pt\vrule width6pt height7pt
depth1pt  \hskip1pt}\bigskip}

\begin{center}

\vspace*{1 cm}

{\huge Convex ancient solutions of the mean curvature flow}
\\ \mbox{ } \\ \mbox{ } \\
{\sc  Gerhard Huisken} \\
{\small Fachbereich Mathematik, Universit\"at T\"ubingen}
\\ {\small Auf der Morgenstelle 10, 72076 T\"ubingen, Germany} \\
{\small e-mail: gerhard.huisken@uni-tuebingen.de} \\
\mbox{ } \\
{\sc  Carlo Sinestrari}
\\ {\small
Dipartimento di Matematica, Universit\`{a} di Roma
``Tor Vergata''}
\\
{\small Via della Ricerca Scientifica, 00133 Roma, Italy}
\\
{\small e-mail: sinestra@mat.uniroma2.it}
\end{center}
\date{}

\vspace{1 cm}

\noindent{\bf Abstract:} We study solutions of the mean curvature flow which are defined for all negative curvature times, usually called ancient solutions. We give various conditions ensuring that a closed convex ancient solution is a shrinking sphere. Examples of such conditions are: a uniform pinching condition on the curvatures, a suitable growth bound on the diameter or a reverse isoperimetric inequality. We also study the behaviour of uniformly $k$-convex solutions, and consider generalizations to ancient solutions immersed in a sphere.

\section{Introduction}\hspace{5 mm}

In this paper, we study convex ancient solutions to the mean curvature flow. We recall that a family of smooth immersions $F:\mm \times (t_0,t_1) \to \rl^{n+1}$, where $\mm$ is an $n$-dimensional manifold, is called a mean curvature flow if it satisfies
 \begin{equation}\label{1.1}
\frac{\partial F}{\partial t}(p,t)=-H(p,t) \nu(p,t),
\qquad p \in \mm, t \in (t_0,t_1),
\end{equation}
where $H(p,t)$ and $\nu(p,t)$ are the mean curvature and the
outer normal respectively at the point $F(p,t)$ of the surface
$\mm_t=F(\cdot,t)(\mm)$. The signs are chosen such that $-H
\nu=\vec{H}$ is the mean curvature vector and the mean
curvature of a convex surface is positive. We assume that $n \geq 2$ and that $\mm$ is closed. 

A solution is called {\em ancient} if it is defined for $t \in (-\infty,0)$. While mean curvature flow is parabolic and hence in general ill-posed backward in time, ancient solutions are of interest for several reasons:
They arise as tangent flows near singularities of the flow and therefore model the asymptotic profile of the surfaces as they approach a singularity, see e.g. \cite{HS1}. Notice that the convexity assumption made in the current paper arises naturally since it was shown in \cite{HS2,Wh} that the rescaling of any singularity of a mean-convex solution of mean curvature flow is (weakly) convex. Ancient solutions of mean curvature flow have also been of interest in theoretical physics where they appear as steady state solutions of boundary renormalisation-group-flow in the boundary sigma model \cite{BS, LVZ}.

Examples of ancient solutions include all homothetically shrinking solutions, in particular the shrinking round sphere $\mm_t = S^n_{R(t)}, R(t) = \sqrt{-2nt}$ and the shrinking cylinders $\mm_t = S^{n-k}_{R(t)}\times \rl^k, R(t) = \sqrt{-2(n-k)t}$. An important example of a non-homothetically shrinking ancient solution is the {\em Angenent oval} \cite{An}, an ancient convex solution of curve shortening in the plane which arises by gluing together near $t \to -\infty$  two opposite translating (non-compact) solutions of curve shortening flow in the plane given by
$$
y_1(t) = - \log\cos x + t, \quad y_2(t) =  \log\cos x - t. 
$$
The translating solution above is known as the {\em grim reaper} curve in the mathematical community while it is known as the {\em hairpin} solution in the physics community, \cite{BS}. The {\em Angenent oval} is known as the {\em paperclip} solution in the physics literature \cite{LVZ}. Compact convex ancient solutions in the plane have been completely classified by Daskalopoulos-Hamilton-Sesum \cite{DHS1} to be either shrinking round circles or an Angenent oval.

A solution analogous to the Angenent oval was constructed in higher dimensions by White \cite{Wh}, p.134. Haslhofer and Hershkovitz give a more detailed construction \cite{HH} and formal asymptotics of this solution was studied by Angenent in \cite{An2}. Heuristic arguments suggest that in higher dimensions many more compact convex ancient solutions may be constructed by appropriately gluing together lower dimensional translating solutions for $t \to -\infty$. A class of such convex eternal solutions was studied by X.-J.-Wang in\cite{Wa}.

We note that there are analogous phenomena for ancient solutions of Ricci flow. A classification of 2-dimensional ancient solutions of positive Gauss curvature was obtained by Daskalopoulos-Hamilton-Sesum in \cite{DHS2}. In higher dimensional Ricci flow Brendle-Huisken-Sinestrari \cite{BHS} proved that the shrinking sphere solution is rigid in a natural class of positively curved Riemannian metrics satisfying a pinching condition.

The current paper first proves a new interior estimate in time for ancient compact convex solutions in theorem 3.1. It then gives a detailed list of characterisations of the shrinking sphere in terms of several natural geometric conditions. These include pinching of the second fundamental form, a scaling invariant diameter bound and a reverse isoperimetric inequality. We summarise these results in section 3 and 4 as follows:

\begin{Theorem}\label{summary}
Let $\mm_t$ be a closed convex ancient solution of mean curvature flow. Then the following properties are equivalent:
\begin{description}
\item[(i)] $\mm_t$ is a family of shrinking spheres.
\item[(ii)] The second fundamental form of $\mm_t$ satisfies the pinching condition $h_{ij} \geq \ep H g_{ij}$ for some $\ep>0$.
\item[(iii)] The diameter of $\mm_t$ satisfies
${\rm diam}(\mm_t) \leq C_1(1+ \sqrt{-t})$ for some $C_1>0$.
\item[(iv)] The outer and inner radius of $\mm_t$ satisfy $\rho_+(t) \leq C_2 \rho_-(t)$ for some $C_2>0$.
\item[(v)] $\mm_t$ satisfies $\max H(\cdot,t) \leq C_3  \min H(\cdot,t)$ for some $C_3>0$.
\item[(vi)]  $\mm_t$ satisfies the reverse isoperimetric inequality $|\mm_t|^{n+1} \leq C_4 |\Omega_t|^{n}$ for some $C_4>0$, where $\Omega_t$ is the region enclosed by $\mm_t$.
\item[(vii)] $\mm_t$ is of type I, that is, $\limsup_{t \to -\infty} \sqrt{-t} \max H(\cdot,t) < \infty$.
\end{description}
\end{Theorem}

Haslhofer and Hershkovitz \cite{HH} have proved a similar result (equivalence of (i), (ii), (iii) and (vii)) under the additional assumption that the solutions are $\alpha$-noncollapsed in the sense of Andrews \cite{A2}. Their method of proof is independent and is based on the recent results by Haslhofer and Kleiner \cite{HK}.

The formal analysis of \cite{An2} supports the existence of ancient solutions whose diameter grows with rate $\sqrt{|t| \ln |t|}$. This suggests that the growth rate in assumption (iii) is close to being optimal.

In section 5 we study ancient convex solutions which are in addition uniformly $k$-convex, i.e. the sum of the smallest $k$ principal curvatures is everywhere bounded below by a fixed fraction of the mean curvature. Such solutions arise when studying surgery algorithms for $k$-convex surfaces, see \cite{HS3}. We prove in section 5 that such solutions satisfy a sharp upper bound on the norm of their second fundamental form establishing ``gaps'' in the geometric behaviour of ancient solutions depending on the level of $k$-convexity. Finally, in section 6 we consider ancient solutions immersed in the sphere, showing that the round shrinking solution can be characterized by curvature pinching conditions which are weaker than in the Euclidean case.

\section{Preliminaries}
\setcounter{equation}{0}

If $\Omega \subset \rl^{n+1}$ is a compact set with nonempty interior, the {\em inner radius} and {\em outer radius} of $\Omega$ are defined as
$$
\rho_-(\Omega)=\max \{ r>0 ~:~ \exists B_r \mbox{ ball of radius $r$} \mbox{ s.t. } B_r \subset \Omega \}, 
$$
$$
\rho_+(\Omega)=\min \{ R>0 ~:~ \exists B_R\mbox{ ball of radius $R$} \mbox{ s.t. } \Omega \subset \overline B_R \}.
$$
In addition, given a unit vector $\nu \in \rl^{n+1}$, the quantity
$$
w(\nu,\Omega)=\max \{ \langle y-x, \nu \rangle ~:~ x,y \in \Omega  \}
$$
is called the {\em width} of $\Omega$ in the direction of $\nu$. By definition, $w(\nu,\Omega)$ is the distance between the two hyperplanes orthogonal to $\nu$ touching $\Omega$ from outside. We then define
$$
w_-(\Omega) = \min_{|\nu|=1} w(\nu,\Omega), \qquad w_+(\Omega) = \max_{|\nu|=1} w(\nu,\Omega).
$$
It is easy to see that
\begin{equation}\label{edw}
w_+(\Omega)= \mbox{diam}(\Omega).
\end{equation}
If $\Omega \subset \rl^{n+1}$ is convex we have the following inequalities \cite[Lemma 5.4]{A}:
\begin{equation}\label{width}
\rho_+ (\Omega) \leq \frac{w_+ (\Omega)}{ \sqrt 2}, \qquad
\rho_-  (\Omega)\geq \frac{w_-(\Omega)}{n+2}.
\end{equation}

If $\mm \subset \rl^{n+1}$ is a smooth closed embedded $n$-dimensional hypersurface, we consider the compact  set $\Omega$ such that $\mm = \partial \Omega$ and define $\rho_\pm(\mm)$ and $w_\pm(\mm)$ to be equal to the corresponding quantities associated to $\Omega$. In addition, we denote by $\mbox{diam}_I (\mm)$ the {\em intrinsic} diameter of $\mm$, computed using the riemannian distance on $\mm$ induced by the immersion, in contrast with the {\em extrinsic} diameter $\mbox{diam} (\mm)$, which is defined in terms of the distance in $\rl^{n+1}$. Elementary arguments, together with \eqref{width}, show that if $\mm$ is convex then
\begin{equation}\label{eidor}
\sqrt{2} \, \rho_+(\mm) \leq  {\rm diam}(\mm) \leq {\rm diam}_I(\mm) \leq \pi \rho_+(\mm).
\end{equation}

For a smooth embedded hypersurface $\mm \subset \rl^{n+1}$, we denote the metric by $g =
\{g_{ij}\}$, the surface measure by $d\mu$, the second fundamental form by $A = \{h_{ij}\}$
and the Weingarten operator by $W=\{h^i_j\}$. We
then denote by
$\ld_1\leq \dots \leq \ld_n$ the principal curvatures, i.e. the eigenvalues
of $W$, and by $H=\ld_1+\dots+\ld_n$ the mean curvature. In addition,
$|A|^2=\ld_1^2+\dots+\ld_n^2$ denotes the squared norm of $A$.

If $\mm_t$ is a family of hypersurfaces evolving by mean curvature flow \eqref{1.1}, the above quantities depend on $t$ and
satisfy the following evolution equations computed in \cite{Hu1}:
\begin{equation}\label{eq2a}
\frac{\partial}{\partial t} h^i_j = \Delta h^i_j + |A|^2 h^i_j,
\end{equation}
\begin{equation}\label{eq2}
\frac{\partial H}{\partial t} = \Delta H + |A|^2 H,
\end{equation}
\begin{equation}\label{eq2b}
 \frac{\partial}{ \partial t} |A|^2 = \Delta |A|^2 - 2|\nabla A|^2 +2|A|^4,
\end{equation}
\begin{equation}\label{eq2c}
\frac{\partial }{\partial t} \, d\mu=-H^2 d\mu. \quad 
\end{equation}

Throughout the paper the surfaces $\mm_t$ have dimension $n \geq 2$, are compact, convex and defined for $t \in (-\infty,0)$. Moreover, $0$ is assumed to be the singular time of the flow, and the surfaces $\mm_t$ shrink to a point as $t \to 0$ by the results of \cite{Hu1}.

The strong maximum principle for tensors applied to equation \eqref{eq2a} implies that all principal curvatures are strictly positive everywhere. In addition, if we consider the evolution equation for the mean curvature \eqref{eq2} together with the inequalities
\begin{equation}\label{eq1}
\frac{H^2}{n} \leq |A|^2 \leq H^2,
\end{equation}
we obtain, by a standard comparison argument,
\begin{equation}\label{eq3}
\min_{\mm_t} H \leq \frac{\sqrt n}{\sqrt{-2t}}, \qquad
\max_{\mm_t} H \geq \frac{1}{\sqrt{-2t}}, \qquad \forall \, t \in (-\infty,0).
\end{equation}

Comparison with evolving spheres, together with the property that $\rho_-(\mm_t) \to 0$ as $t \to 0$, yields easily the following bounds on the innner and outer radius of $\mm_t$:
\begin{equation}\label{eq4}
\rho_-(\mm_t) \leq \sqrt{-2nt} \leq \rho_+(\mm_t), \qquad \forall \, t \in (-\infty,0).
\end{equation}

We finally recall Hamilton's differential Harnack estimate \cite{Ha3}, which for ancient solutions takes the form
\begin{equation}\label{e.HH}
\frac{\partial H}{\partial t} - \frac{|\nabla H|^2}{H} \geq 0.
\end{equation}
In particular, this implies that $H$ is pointwise nondecreasing. Therefore, our solutions have uniformly bounded curvature on any time interval of the form $(-\infty,T_1]$, with $T_1<0$. Since $H$ is the speed of our evolving surfaces, we deduce that for each solution there exists a constant $K >0$, such that the outer radius satisfies the bound 
\begin{equation}\label{ebh}
\rho_+(\mm_t) \leq K(1+|t|) \qquad \forall \, t<0.
\end{equation}

\section{Pinched solutions}
\setcounter{equation}{0}

In this section we consider ancient solutions satisfying the pinching condition
\begin{equation}\label{e.hyp}
h_{ij} \geq \ep H g_{ij}
\end{equation}
for some $\ep>0$ independent of $t$. We will show that a solution with this property is necessarily a family of shrinking spheres.
To this purpose we consider the function, introduced in \cite{Hu1}, 
\begin{equation}\label{gsigma}
f_{\sigma}=\frac{|A|^2-H^2/n}{H^{2-\sigma}},
\end{equation}
with $\sigma >0$. The function $f_\sigma$ is nonnegative, and it is zero only at umbilical points. We first prove an integral estimate for this function which holds on any pinched solution of the flow, not necessarily ancient, and only depends on the pinching constant and the lifespan of the solution. The characterization of ancient solutions will follow as a corollary.

\begin{Theorem}\label{t01}
Let $\mm_t$, with $t \in [T_0,0)$ be a solution to the mean curvature flow such that $\mm_{T_0}$ is closed, convex and satisfies \eqref{e.hyp} for some $\ep>0$, and which becomes singular as $t \to 0^-$.Then there exist $c_1,c_2,c_3>0$ depending only on $n,\ep$ such that, for every $p,\sigma>0$ satisfying
$$p \geq c_1, \qquad \sigma \leq c_2/\sqrt p, \qquad p \sigma >n,$$
we have that
$$
\left( \int_{\mm_t} f_\sigma^p \, dt \right) ^\frac{2}{\sigma p} \leq  \frac{c_3}{|T_0|^{1- \frac{n}{\sigma p}} - |t|^{1- \frac{n}{\sigma p}}},
$$
for all $t \in [T_0,0)$.
\end{Theorem}
{\bf Proof} --- As shown in \cite{Hu1}, the solution $\mm_t$ is convex and satisfies \eqref{e.hyp} for all $t \in (T_0,0)$. In addition, by an easy adaptation of the proof of Lemma 5.5 of \cite {Hu1} one finds that, if $p,\sigma$ are such that
$$
p \geq \frac{100}{\ep^2}, \qquad 
\sigma \leq \frac{n \ep^3}{16} \frac{1}{\sqrt p}
$$
then
\begin{equation}\label{e1.1}
{\df{d}{dt}} \int_{\mm_t} f_\sigma^p d\mu  \leq  - p\sigma
\int_{\mm_t} H^2 f_\sigma^p d\mu.
\end{equation}
Recalling inequalities \eqref{eq1} we see that $0 \leq f < H^\sigma$ and therefore
\begin{equation}\label{e1.2}
{\df{d}{dt}} \int_{\mm_t} f_\sigma^p d\mu  \leq  - p\sigma
\int_{\mm_t} f_\sigma^{p+\frac 2\sigma} d\mu \leq - p\sigma \left( \int_{\mm_t} f_\sigma^p d\mu \right)^{1+\frac 2{\sigma p}} |\mm_t|^{-\frac{2}{\sigma p}},
\end{equation}
by  H\"older's inequality, where $|\mm_t|$ denotes the area of $\mm_t$. To exploit this estimate, we need an upper bound on $|\mm_t|$. We first observe that, by property (\ref{e.hyp}), the Gauss curvature $K$ of our surfaces satisfies $K \geq (\ep H)^n$. Therefore, by the Gauss-Bonnet theorem
$$
\int_{\mm_t} H^n \, d\mu \leq \ep^{-n} \int_{\mm_t} K \, d\mu = c(\ep,n).
$$
Therefore, by \eqref{eq2c},
$$
\frac {d}{dt} |\mm_t| = -\int_{\mm_t} H^2 \, d\mu \geq -\left(\int_{\mm_t} H^n d\mu \right)^{2/n} |\mm_t|^{1-\frac 2n} \geq -c|\mm_t|^{1-\frac 2n},
$$
where we denote by $c$ a suitable positive constant depending on $\ep,n$ (possibly different from one formula to another). By integrating this inequality over $[t,s]$, with $T_0 \leq t<s<0$ arbitrary, we obtain
$$
|\mm_s|^{2/n} - |\mm_t|^{2/n} \geq - c (s-t).
$$
We have $|\mm_s| \to 0$ as  $s \to 0$ by the result of \cite{Hu1}. Therefore, letting $s \to 0$ we obtain
\begin{equation}\label{e1.3}
|\mm_t| \leq c (-t)^{n/2}.
\end{equation}
Let us set for simplicity $\phi(t):=\int_{\mm_t} f_\sigma^p \, d\mu$. Combining (\ref{e1.2}) and (\ref{e1.3}) we obtain
\begin{equation}\label{eq.av}
{\df{d}{dt}}  \phi^{-\frac{2}{\sigma p}} \geq c (-t)^{- \frac{n}{\sigma p}},
\end{equation}
for any $t$ such that $\phi(t)>0$. Observe that $\phi(t) \geq 0$, with equality if and only if $\mm_t$ is a sphere. This implies in particular that, if $\phi(t) > 0$, then $\phi(s) >0$ for all $s <t$. Therefore, if we take any $t \in (T_0,0)$ such that $\phi(t)>0$ we can integrate \eqref{eq.av} over $[T_0,t]$ and obtain
\begin{eqnarray*}
\phi^{-\frac{2}{\sigma p}} (t) & \geq & \phi^{-\frac{2}{\sigma p}} (T_0) + c \int_{|t|}^{|T_0|} \tau^{- \frac{n}{\sigma p}} \, d \tau \\
& > &  c \int_{|t|}^{|T_0|} \tau^{- \frac{n}{\sigma p}} \, d \tau > c \left( |T_0|^{1-\frac{n}{\sigma p}}-|t|^{1-\frac{n}{\sigma p}} \right),
\end{eqnarray*}
where we have also used the assumption $\sigma p > n$. This proves our assertion if $\phi(t)>0$, while the case $\phi(t)=0$ is trivial. \qed

As a corollary, we obtain the following result which gives the equivalence between (i) and (ii) in Theorem \ref{summary}.

\begin{Theorem}\label{t1}
Let $\mm_t$, with $t \in (-\infty,0)$ be an ancient solution to the mean curvature flow such that every $\mm_t$ is closed, convex and satisfies \eqref{e.hyp} for some $\ep>0$ independent of $t$. Then $\mm_t$ is a family of shrinking spheres.
\end{Theorem}
{\bf Proof} --- We can let $T_0 \to -\infty$ in the previous theorem and conclude that for suitable values of $\sigma,p$ the integral $\int_{\mm_t} f_\sigma^p $ is zero for every $t<0$. This implies that every $\mm_t$ is a sphere. \qed

\section{Solutions with a diameter bound}
\setcounter{equation}{0}

In this section we prove the remaining parts of Theorem \ref{summary}. We first show that a growth bound of the order $O(\sqrt{-t})$ on the diameter of the solution allows to control the variation of the curvature at any fixed time.

\begin{Lemma}\label{l2}
Let $\mm_t$, with $t \in (-\infty,0)$, be a closed convex ancient solution of the mean curvature flow. Then the two following properties are equivalent.
\begin{description}
\item[(i)] There exists $C_1>0$ such that
\begin{equation}\label{hyp2}
{\rm diam} (\mm_t) \leq C_1(1+ \sqrt{-t}), \quad \mbox{ for all }t<0.
\end{equation}
\item[(ii)]
There are $C',C''>0$ such that
\begin{equation}\label{hyp4}
\frac{C'}{\sqrt{-t}} \leq H \leq \frac{C''}{\sqrt{-t}},  \quad \mbox{ on } \mm_t, \mbox{ for all }t<0.
\end{equation}
\end{description}
\end{Lemma}
{\bf Proof} --- Let us show that \eqref{hyp2} implies \eqref{hyp4}. We first observe that, by property \eqref{eidor}, the intrinsic diameter of our surfaces satisfies
\begin{equation}\label{hyp2a}
{\rm diam}_I (\mm_t) \leq C_1\frac{\pi}{\sqrt 2} \left(1+ \sqrt{-t} \, \right), \ \mbox{ for all }t<0.
\end{equation}
Then  we also have 
\begin{equation}\label{hyp3}
{\rm diam}_I(\mm_t) \leq c \, \sqrt{-t}, \ \mbox{ for all }t<0,
\end{equation}
for a suitable $c>0$. In fact, for $t$ close to zero this follows from the convergence of $\mm_t$ to a round point \cite{Hu1}, while for $t$ away from zero it follows from \eqref{hyp2a}.

Next we recall a well known consequence of Hamilton's Harnack inequality \eqref{e.HH}, stating that
\begin{equation}\label{e3.10}
H(p_1,t_1) \leq H(p_2,t_2) \, \exp \left( \frac{{\rm diam}_I^2(M_{t_1})}{4(t_2-t_1)} \right)
\end{equation}
for any $p_1,p_2 \in \mm$ and $t_1<t_2<0$. For any $t<0$, we apply this inequality, together with \eqref{hyp3}, with $t_1=t$, $t_2=t/2$, and we obtain
$$
\max H(\cdot ,t) \leq  e^{c^2/2} \min H \left( \cdot ,  t/2 \right).
$$
Using (\ref{eq3}), we find
$$
\max H(\cdot ,t) \leq e^{c^2/2}  \min H \left( \cdot ,  t/2 \right) \leq e^{c^2/2} \sqrt\frac{n}{-t}, \qquad \forall \, t <0.
$$
Since $t<0$ is arbitrary, we can replace $t$ by $2t$ and use the other part of \eqref{eq3} to obtain
$$
\min H \left( \cdot ,  t \right) \geq e^{-c^2/2}\max H(\cdot ,2t) \geq e^{-c^2/2} \frac{1}{2\sqrt{-t}}, \qquad \forall \,  t <0.
$$
The two inequalities together imply \eqref{hyp4}. 

Suppose now that \eqref{hyp4} holds. Since $\mm_t$ shrinks to a point as $t \to 0$, we find for any pair of points $p,q \in \mm$
$$
|F(p,t)-F(q,t)| \leq \int_t^0 H(p,\tau) \, d\tau + \int_t^0 H(q,\tau) \, d\tau \leq 2 C'' \int_t^0 \frac{d\tau}{\sqrt{-\tau}} =  4 C'' \sqrt{-t},
$$
which implies \eqref{hyp2}. \qed

\begin{Theorem}\label{t.db}
Let $\mm_t$, with $t \in (-\infty,0)$, be a closed convex ancient solution of the mean curvature flow satisfying either \eqref{hyp2} or \eqref{hyp4}. Then $\mm_t$ is a family of shrinking spheres.
\end{Theorem}
{\bf Proof} --- If our solution satisfies $h_{ij} \geq \ep H g_{ij}$ for a time-independent $\ep>0$, then the assertion follows from Theorem \ref{t1}. We assume therefore that this does not hold; then there exists a sequence $\{(p_k,t_k)\}$, with $t_k \to -\infty$, such that $(\ld_1/H)(p_k,t_k) \to 0$ as $k \to \infty$. We consider the flow $\mm_t$ for $t \in [2t_k,t_k]$ and rescale it by a factor $1/\sqrt{|t_k|}$ in space and $1/|t_k|$ in time. We obtain a sequence of flows which are all defined for $t \in [-2,-1]$ and, by the previous lemma, have curvature and diameter which are uniformly bounded from above and below by positive constants. By \cite{EH}, we also have uniform bounds on all derivatives of the curvature for $t \in [-3/2,-1]$. So we can extract a subsequence converging smoothly to a solution of the flow in $[-3/2,-1]$. This limit solution is convex and compact, but it contains a point with $\ld_1=0$ at $t=-1$. A well known argument based on the strong maximum principle, see \cite{Ha2}, implies that the limit solution must split as a product containing a flat factor,  in contradiction with the diameter bound. This shows that the original solution $\mm_t$ satisfies $h_{ij} \geq \ep H g_{ij}$, and we can conclude by Theorem \ref{t1}. \qed

\begin{Corollary}\label{c4.3} If our ancient solution $\mm_t$ satisfies either of the two properties:
 $$\rho_+(\mm_t) \leq C \rho_-(\mm_t), \ \mbox{ for all }t<0
 $$
 $$
 \max H(\cdot,t) \leq C \min H(\cdot,t), \ \mbox{ for all }t<0
 $$
 for a constant $C>0$, then $\mm_t$ is a family of shrinking spheres.
\end{Corollary}
{\bf Proof} --- From inequalities (\ref{eidor}), (\ref{eq4}) we see that the hypothesis on $\rho_-,\rho_+$ implies the diameter bound \eqref{hyp2}. Also, the assumption on $H$ implies, by \eqref{eq3}, that \eqref{hyp4} holds. So in both cases we can conclude by
the previous theorem. \qed

We now consider an assumption of a different kind on our ancient solution $\mm_t$. As in Section 2, we denote by $\Omega_t$ the compact region enclosed by $\mm_t$. We say that $\mm_t$ satisfies a {\em uniform reverse isoperimetric estimate} if there exists a constant $C>0$ such that
\begin{equation}\label{rii}
|\mm_t|^{n+1} \leq C |\Omega_t|^n, \qquad \forall t<0,
\end{equation}
where $|\mm_t|$ and $|\Omega_t|$ denote the $n$-dimensional and $(n+1)$-dimensional measure of $\mm_t$ and $\Omega_t$ respectively. Clearly, the constant $C$ in \eqref{rii} has to be greater than the optimal constant in the isoperimetric inequality achieved by the sphere. The next lemma shows that such an assumption implies a uniform bound on the ratio between outer and inner radius.

\begin{Lemma}
For any $n \geq 1$ and $c_1 >0$ there exists $c_2=c_2(c_1,n) \geq 1$ with the following property. Let $\mm \subset \rl^{n+1}$ be a closed convex $n$-dimensional hypersurface such that
$$
|\mm|^{n+1} \leq c_1 |\Omega|^{n},
$$
where $\Omega$ is the region enclosed by $\mm$. Then the outer and inner radius $\rho_+$ and $\rho_-$ of $\mm$ satisfy
$$
\frac{\rho_+}{\rho_-} \leq c_2.
$$
\end{Lemma}
{\bf Proof} --- We shall use the notion of maximal and minimal width $w_+(\mm)$ and $w_-(\mm)$ recalled in Section 2; since $\mm$ is fixed throughout the proof, we simply write $\rho_\pm,w_\pm$ instead of $\rho_\pm(\mm), w_\pm(\mm)$.

Assume for simplicity that the direction achieving the minimal width $w_-$ is the one of the $x_{n+1}$--axis, and denote by $\Sigma$ the orthogonal projection of $\mm$ onto the $\{x_{n+1}=0\}$ hyperplane. Then we can estimate 
$$
|\Omega| < w_- |\Sigma| \qquad |\mm| > |\Sigma|,
$$
where $|\Sigma|$ is the $n$-dimensional measure of $\Sigma$. 
Thus the hypothesis of our lemma implies
\begin{equation}\label{er.1}
|\Sigma| \leq c_1 w_-^n.
\end{equation}
In addition, we observe that if $P,Q$ are any two points in $\mm$ and $P',Q' \in \Sigma$ are their projections onto the $\{x_{n+1}=0\}$ hyperplane, we have
$$
|P-Q| \leq |P'-Q'| + w_-.
$$
Therefore, recalling also \eqref{edw}, 
\begin{equation}\label{er.2}\mbox{diam}(\Sigma) \geq \mbox{diam}(\mm) - w_- = w_+-w_-.
\end{equation}
In the case $n=1$ we have that ${\rm diam}(\Sigma)=|\Sigma|$; hence, \eqref{er.1} and \eqref{er.2} together yield
$w_+ \leq (c_1+1)w_-$, which implies the assertion, by \eqref{width}.

If $n>1$ we argue as follows. We take an $(n+1)$-dimensional ball $B_0$ of radius $\rho_-$ contained in $\Omega$ and let $B_1$ be its projection on the $x_{n+1}=0$ hyperplane; then $B_1 \subset \Sigma$. Let $P_1 \in \Sigma$ be the center of $B_1$. By definition of diameter, there exists $P_2 \in \Sigma$ such that $|P_2-P_1| \geq \mbox{diam}(\Sigma)/2$. We then intersect $B_1$ with the hyperplane through $P_1$ orthogonal to the direction $P_2-P_1$, obtaining the $(n-1)$-dimensional ball
$$
B_2=\{ P \in B_1 ~:~ \langle P-P_1 , P_2 -P_1 \rangle = 0 \}.
$$
Since $\Sigma$ is convex and it contains both $B_2$ and $P_2$, it also contains the cone $K$ with basis $B_2$ and vertex $P_2$. We have that $B_2$ is an $(n-1)$-dimensional ball of radius $\rho_-$, while the height of $K$ equals $|P_2-P_1| \geq \mbox{diam}(\Sigma)/2$, and therefore
$$
|\Sigma| \geq |K| \geq \frac{\omega_{n-1}}{n} \rho_-^{n-1} \frac{\mbox{diam}(\Sigma)}{2},
$$
where $\omega_{n-1}$ is the volume of the unit $(n-1)$-dimensional ball.
Using \eqref{width} and \eqref{er.2} we deduce
$$
|\Sigma| \geq \frac{\omega_{n-1}}{n} \frac{w_-^{n-1}}{(n+2)^{n-1}} \frac{w_+-w_-}{2} = :\kappa_n w_-^{n-1} (w_+-w_-).
$$
where $\kappa_n$ only depends on $n$. Then \eqref{er.1} implies
$$
c_1 w_- \geq \kappa_n (w_+-w_-),
$$
which gives $w_+ \leq (1 + c_1/\kappa_n) w_-$. By \eqref{width}, we obtain the conclusion. \qed

Combining the previous lemma with Corollary \ref{c4.3}, we obtain

\begin{Corollary}
Suppose that there exists a constant $C>0$ such that the uniform reverse isoperimetric estimate \eqref{rii} holds. Then $\mm_t$ is a family of shrinking spheres.
\end{Corollary}

In analogy with the terminology used in the study of finite time singularities, we say that an ancient solution to the mean curvature flow is of type I if there exist constants $C>0$ and $T_0<0$ such that
\begin{equation}\label{typeI}
\max_{\mm_t}H(\cdot,t) \leq C/\sqrt{-t}, \ \mbox{ for all }t\leq T_0.
\end{equation} 
Otherwise we call it of type II. 

\begin{Proposition}
A closed convex ancient solution of the mean curvature flow of type I is a family of shrinking spheres.
\end{Proposition}
{\bf Proof} --- The results of \cite{Hu1} imply that inequality \eqref{typeI} also holds for $t \in [T_0,0)$ possibly with a larger $C$. Then we can argue as in the last step of the proof of Lemma \ref{l2}, to obtain that a type I solution satisfies the bound \eqref{hyp2}. We conclude by Theorem \ref{t.db}. \qed 

The previous result completes the proof of Theorem \ref{summary}. To conclude the section, we show that any ancient solution different from the sphere can be rescaled in order to produce a translating soliton, i.e. a translating solution of the flow. In fact, the type II property allows to perform a procedure analogous to the case finite time singularities of the mean curvature flow, see e.g. \cite{HS1}.

\begin{Theorem}
If $\mm_t$ is a closed convex ancient solution of the flow different from a shrinking sphere, then there is a family of rescaled flows which converge to a translating soliton.
\end{Theorem}
{\bf Proof} --- By the previous proposition, $\mm_t$ is of type II. For any given integer $k>1$, choose $p_k \in \mm$ and $t_k \in [-k,-1]$ such that
$$
\sqrt{-t_k} H(p_k,t_k) =\max_{p \in \mm, t \in [-k,-1]} \sqrt{-t} H(p,t),
$$
and set $L_k=H(p_k,t_k)$. The type II property implies that
$$
t_k \to -\infty, \qquad \sqrt{-t_k} L_k \to +\infty 
$$
as $k \to +\infty$. We remark that, in contrast to the rescaling near a finite time singularity, in this case $L_k$ does not go to $+\infty$ and in general may go to zero. Now we define a rescaled immersion as follows
 \begin{equation}\label{res}
F_k(\cdot,\tau)=L_k
F \left(\ \cdot \ ,\frac{\tau}{L_k^2}+t_k \right), \qquad \tau \in (-\infty,L_k^2(-t_k)),
 \end{equation}
and we denote by $H_k$ the mean curvature on the $k$-th rescaled flow. Observe that 
$L_k^2(-t_k) \to + \infty$ as $k \to \infty$. In addition, we have $H_k(p_k,0)=1=\max H(\cdot,0)$.

With computations similar to \S 4 in \cite{HS1}, we obtain that for any $\ep>0$ small and for any $\overline T>0$ large, the rescaled flows satisfy $H_k \leq 1+\ep$ on the time interval $(-\infty,\overline T)$ for $k$ large enough. This ensures convergence of a subsequence of the rescaled flows to a limit solution defined for all times, convex, with $H \leq 1$ everywhere and $H =1$ at some point for $t=0$. By a result of Hamilton \cite{Ha3}, it must be a translating soliton. \qed

\section{Uniformly k-convex ancient solutions}
\setcounter{equation}{0}

We recall that a hypersurface is called $k$-convex, for some $k \in \{1,\dots,n-1\}$, if $\ld_1+\dots+\ld_k \geq 0$ everywhere. Thus, we say that a solution of the mean curvature flow is uniformly $k$-convex if it satisfies
\begin{equation}\label{hypk}
\ld_1+\dots+\ld_k \geq \alpha H >0
\end{equation}
for some $\alpha>0$ and for all $t<0$. Such an inequality is preserved by the flow. Observe that uniform $1$-convexity means convexity with curvature pinching, as in assumption \eqref{e.hyp}. For general $k$, we can give a sufficient condition for $k$-convexity in terms of the quotient $|A|^2/H^2$.

\begin{Lemma}\label{lkc}
Suppose that at a point of an $n$-dimensional hypersurface $\mm$ we have $H>0$ and 
\begin{equation}\label{condkconv}
\frac{|A|^2}{H^2} \leq \frac{1-2\alpha}{n-k}
\end{equation}
for some $k=1,\dots,n-1$ and some $\alpha \in (0,1/2)$. Then at the same point we also have
$$
\ld_1+\dots+\ld_k \geq \alpha H.
$$
\end{Lemma}
{\bf Proof} --- 
Let us set
$
\mu=\ld_{k+1}+\dots+\ld_n.
$
Then we have
$$\mu^2 \leq (n-k) (\ld_{k+1}^2+\dots+\ld_n^2) \leq (n-k) |A|^2.
$$
Therefore, by \eqref{condkconv},
\begin{eqnarray*}
(1-\alpha)^2 H^2   >  (1-2\alpha)H^2 \geq  (n-k)|A|^2 \geq \mu^2.
\end{eqnarray*}
We conclude 
$$
\ld_1+\dots+\ld_k =H-\mu \geq H - (1-\alpha)H = \alpha H. \qed
$$

In \cite{HS3} we proved that a $2$-convex solution (not necessarily convex) of the mean curvature flow on a finite time interval $[0,T)$ satisfies the following property: for any $\eta>0$ there exists $C_\eta>0$ such that
$$
|A|^2-\frac{1}{n-1}H^2 \leq \eta H^2 + C_\eta.
$$
The proof can be easily generalized to show that $k$-convex solutions satisfy an analogous inequality with $\frac{1}{n-1}$ replaced by $\frac{1}{n-k+1}$. Here we show that an ancient solution which is convex and uniformly $k$-convex satisfies a stronger estimate without the terms in the right hand side.

\begin{Theorem} \label{tkc}
Let $\mm_t$, with $t \in (-\infty,0)$, be a convex closed ancient solution of the mean curvature flow embedded in $\rl^{n+1}$ with $n \geq 3$, satisfying \eqref{hypk} for some $k=2,\dots,n-1$ and $\alpha>0$.
Then we have $H^2 > (n-k+1) |A|^2$ everywhere.
\end{Theorem}
 {\bf Proof} ---In the first part of the proof we follow a procedure similar to section 5 in \cite{HS3}. 
We introduce, for $\eta>0$ and $\sigma \in [0,2]\,$, the function
\begin{equation}\label{gseta}
f_{\sigma,\eta}=\frac{|A|^2-\left(\frac{1}{n-k+1}+\eta \right)H^2}{H^{2-\sigma}}.
\end{equation}
Then
\begin{eqnarray}
\frac{\partial f_{\sigma,\eta}}{\partial t} & = &
\Delta f_{\sigma,\eta} + \frac{2(1-\sigma)}{H}\langle \nb H, \nb
f_{\sigma,\eta} \rangle - {\displaystyle \frac{\sigma(1-\sigma)}
{H^2}f_{\sigma,\eta}}|\nb H|^2 \label{eqg}\\
& & - {\displaystyle \frac{2}{H^{4-\sigma}}}|H\nb_i h_{kl}-\nb_i H \,
h_{kl}|^2 + \sigma|A|^2 f_{\sigma,\eta}. \nonumber
 \end{eqnarray}
In particular, for $\sigma=0$ the maximum of $f_{\sigma,\eta}$ decreases in time, which shows that the inequality $|A|^2 \leq (\frac{1}{n-k+1}+\eta)H^2$ is preserved by the flow.

For simplicity, we write $f \equiv f_{\sigma,\eta}$, and we denote by $f_+$ the positive part of $f$.
From \eqref{eqg} we obtain, as in \cite[Lemma 5.4]{HS3}, that
there exist constants $c_1, c_2>1$, depending only on $n,k$ and the constant $\alpha$ in \eqref{hypk},
such that
\begin{eqnarray}
{\df{d}{dt}} \int_{\mm_t} f_+^p d\mu & \leq & - {\df{p(p-1)}{2}}
\int_{\mm_t} f_+^{p-2} |\nb f|^2 d\mu- {\df{p}{c_1}} \int_{\mm_t}
{\df{f_+^p}{H^2}} |\nb H|^2 d\mu \nonumber \\
& & + p\sigma
\int_{\mm_t} |A|^2 f_+^p d\mu \label{e3.1}
\end{eqnarray}
for any $p \geq c_2$. We now observe that
$$(n-k+1)|A|^2- H^2  -  \sum_{k \leq i<j \leq n} (\ld_i-\ld_j)^2
 =  (n-k) \sum_{i=1}^{k-1} \ld_i^2 -  2  \sum_{i=1}^{k-1} \sum_{j=i+1}^n\ld_i\ld_j. 
 $$
Since for any $i<k$ 
$$
\sum_{j=i+1}^n \ld_j \geq (n-i) \ld_i > (n-k) \ld_i,
$$
we deduce that the right-hand side in the above formula is negative, that is
\begin{equation}
(n-k+1)|A|^2- H^2  <  \sum_{k \leq i<j} (\ld_i-\ld_j)^2.
\end{equation}
Therefore, if we set $Z=H \rm{tr}(A^3)-|A|^4$, we have
\begin{eqnarray*}
Z & = & \sum_{i<j} \ld_i\ld_j(\ld_i-\ld_j)^2 > 
\sum_{k \leq i<j}\ld_i\ld_j(\ld_i-\ld_j)^2  \\
& \geq & 
\sum_{k \leq i<j}\ld_k^2(\ld_i-\ld_j)^2 \geq
\frac{(\alpha H)^2}{k^2}\sum_{k<i<j}(\ld_i-\ld_j)^2 \\
& \geq & \frac{(\alpha H)^2}  {k^2} ((n-k+1) |A|^2- H^2 ).
\end{eqnarray*}
Keeping into account the definition of $f=f_{\sigma,\eta}$, we see that
at the points where $f>0$ we have
\begin{equation}\label{e.zeta}
Z \geq \frac{(n-k+1)\alpha^2 \eta }{k^2} H^4.
\end{equation}
On the other hand, we have (see the proof of Lemma 5.5 in \cite{HS3})
\begin{equation}\label{e3.2}
\int \df{2Z}{H^2}f_+^p d\mu \leq
8np\int \df{f_+^{p-1}}{H}|\nb H|\,|\nb f| d\mu
+10n \int \df{f_+^p}{H^2} |\nabla H|^2 d\mu.
\end{equation}
From estimates \eqref{e3.1}, \eqref{e.zeta} and \eqref{e3.2} we conclude, as in Proposition 5.6 of \cite{HS3}, that there exist constants $c_3,c_4$ such that
\begin{equation}\label{e3.3}
{\df{d}{dt}} \int_{\mm_t} f_+^p d\mu  \leq  - p\sigma
\int_{\mm_t} |A|^2 f_+^p d\mu 
\end{equation}
for all $p \geq c_3$ and $\sigma \leq \eta(c_4 \sqrt p)^{-1}$.

Now let us denote by $\Omega_t$ the region enclosed by $\mm_t$. The volume $|\Omega_t|$ satisfies
$$
\frac{d}{dt}|\Omega_t|= - \int_{\mm_t} H \, d\mu.
$$
From the bound \eqref{ebh} 
we deduce that there exists $c_1>0$ such that
\begin{equation}\label{e.vol}
\int_{t}^0 \left( \int_{\mm_\tau} H \, d\mu \right) \, d\tau =|\Omega_{t}| \leq c_1 (1+|t|)^{n+1}, \qquad \forall \, t<0.
\end{equation}
For any fixed $\eta>0$, let us now choose $p,\sigma$ such that \eqref{e3.3} holds and in addition
\begin{equation}\label{eps}
p \sigma > 2n+1.
\end{equation}
We also set
$$
\gamma=\frac{2}{p\sigma+1}.
$$
By \eqref{eps} we have $\gamma < 1/(n+1)<1$. Let $f_+$ denote the positive part of $f$. Keeping into account that $f \leq H^\sigma$, we find
\begin{eqnarray*}
\int f_+^p \, d\mu &  \leq & \int f_+^{p(1- \gamma)} H^{p\sigma\gamma} \, d\mu \\
& = &  \int (f_+^p H^2 )^{1-\gamma} H^{p \sigma \gamma-2(1-\gamma)} \, d\mu =
\int (f_+^p H^2 )^{1-\gamma} H^\gamma \, d\mu \\
& \leq & \left( \int f_+^p H^2 \,  d\mu\right)^{1-\gamma} \left( \int H \, d\mu \right)^\gamma.
\end{eqnarray*}
Therefore, by \eqref{e3.3} and \eqref{eps},
\begin{eqnarray*}
\df d{dt} \int f_+^p d\mu & \leq & - p \sigma \int |A|^2 f_+^p d \mu \leq -(2n+1)\int |A|^2 f_+^p d \mu   < - \int H^2 f_+^p d \mu \\
& \leq & -  \left( \int f_+^p \,  d\mu\right)^\frac 1{1-\gamma} \left( \int H \, d\mu \right)^{-\frac \gamma{1-\gamma}} \\
& = & - \left( \int f_+^p \,  d\mu\right)^{1+\frac 2{p\sigma-1}} \left( \int H \, d\mu \right)^{-\frac 2{p\sigma-1}}.
\end{eqnarray*}
Let us set for simplicity
$$
\phi(t)=\int f_+^p d\mu, \qquad \psi(t)=  \int H \, d\mu, \qquad \beta=\frac{2}{p\sigma-1}.
$$
Then the previous inequality implies
$$
\frac {d}{dt}(\phi^{-\beta}) = -\beta \phi^{-\beta-1}\frac {d}{dt} \phi \geq \beta \psi^{-\beta}. 
$$
Suppose that there exists $t_1<0$ such that $\phi(t) \neq 0$ for all $t \leq t_1$. Then we have, for any $t_0<t_1$,
$$
\frac{1}{\phi^\beta} (t_1) - \frac{1}{\phi^\beta} (t_0) \geq \beta \int_{t_0}^{t_1} \psi^{-\beta} \, dt \geq \beta (t_1-t_0)^{1+\beta} \left( \int_{t_0}^{t_1} \psi \, dt \right)^{-\beta}.
$$
Recalling (\ref{e.vol}) we find
\begin{equation}\label{ineq}
\frac{1}{\phi^\beta(t_1)} \geq c_2(t_1-t_0)^{1+\beta}(1-t_0)^{-\beta(n+1)}, \quad \forall t_0<t_1 <0,
\end{equation}
for some $c_2>0$. Assumption \eqref{eps} implies that $1-n\beta>0$, and so the right-hand side of (\ref{ineq}) diverges as $t_0 \to -\infty$ giving a contradiction. It follows that we can find negative values of $t$ with $|t|$ arbitrarily large such that $\phi(t)=0$. But $\phi(t)=0$ if and only if $|A|^2 \leq (\frac{1}{n-k+1}+\eta)H^2$ everywhere on $\mm_t$. Such an inequality is preserved by the flow, so it must hold for all $t<0$. By the arbitrariness of $\eta>0$, we obtain that $|A|^2 \leq \frac{1}{n-k+1}H^2$ for every $t < 0$. By the strong maximum principle applied to \eqref{eqg}, we either have $|A|^2 < \frac{1}{n-k+1}H^2$ everywhere, or $|A|^2 \equiv \frac{1}{n-k+1}H^2$ everywhere, which implies that $\mm_t$ is a cylinder; the latter is excluded, since our solution is compact. 
\qed

\begin{Corollary}
Let $\mm_t$, with $t \in (-\infty,0)$, be an ancient solution of the mean curvature flow such that $\mm_t$ is closed and convex for all $t$. Then there is an integer $h=1,\dots,n$ such that
$$
\sup_{t \in (-\infty,0)} \max_{\mm_t} \frac{|A|^2}{H^2} = \frac 1h.
$$
\end{Corollary}
{\bf Proof} --- Since on our solution $|A|^2 \leq H^2 \leq n |A|^2$, the supremum in the statement lies between $1/n$ and $1$. Let us assume that it is smaller than $1$, since otherwise the statement holds. Then there exists $h=1,\dots,n-1$ such that
$$
\sup_{t \in (-\infty,0)} \max_{\mm_t} \frac{|A|^2}{H^2} \in \left[ \frac 1{h+1} , \frac 1h \right).
$$
We can now use Lemma \ref{lkc} to deduce that our surfaces are uniformly $k$-convex, with $k=n-h$. If $k=1$, then our surfaces are spheres by Theorem \ref{t1}, which means that $|A|^2/H^2 \equiv 1/n$. If $k=2,\dots,n-1$, then Theorem \ref{tkc} implies that the supremum of $|A|^2/H^2$ is not greater than $1/(n-k+1) = 1/(h+1)$; by the choice of $h$, it must coincide with $1/(h+1)$. \qed

Under an additional assumption on the rate of curvature decay, we can show that $k$-convex ancient solutions satisfy a gradient estimate similar to Theorem 6.1 in \cite{HS3}, without the constant remainder term, provided $k$ is small enough compared to $n$.

 \begin{Theorem}
Let $\mm_t \subset \rl^{n+1}$, with $t \in (-\infty,0)$, be an ancient solution of the mean curvature flow such that the $\mm_t$'s are closed, convex and uniformly $k$-convex for some $k<(2n+1)/3$. Suppose in addition that 
\begin{equation}\label{hyp4b}
\int_{-\infty}^{-1} \left( \min_{\mm_t} |A|^2 \right) dt=+\infty.
\end{equation}
Then there exists $C=C(k,n)$ such that 
\begin{equation}\label{grest}
|\nabla A|^2 \leq C |A|^4
\end{equation}
on $\mm_t$, for all $t<0$.
\end{Theorem}
{\bf Proof} ---
Let us define
$\sigma=\sigma(n,k)$ as
\begin{equation}\label{kappan}
\sigma=\frac 12 \left(\frac{3}{n+2} -\frac{1}{n-k+1} \right).
\end{equation}
Then $\sigma>0$ by our hypothesis on $k$. We then set
$$
g_1=\left(\frac{1}{n-k+1}+\sigma \right)H^2-|A|^2, \qquad g_2= \frac{3}{n+2} H^2 - |A|^2.
$$
By Theorem \ref{tkc}, we have
$$g_2 > g_1 > \sigma H^2 >0.$$
In particular, the quotient $|\nabla A|^2/g_1g_2$ is well defined.
As in the proof of Theorem 6.1 in \cite{HS3}, we find
\begin{eqnarray}
\lefteqn{\frac{\partial}{\partial t} \left( \frac{|\nabla A|^2}{g_1 g_2}
\right) - \Delta \left(\frac{|\nabla A|^2}{g_1 g_2} \right)
-\frac{2}{g_2} \left\langle \nabla g_2, \nabla \left( \frac{|\nabla A|^2}{g_1 g_2} \right) \right\rangle}
\nonumber
\\ & \leq &
\frac{|\nabla A|^2 \, |A|^2}{g_1 g_2} \left( c_n - 2
\sigma^2\frac{n+2}{3n}
\frac{|\nabla A|^2}{g_1 g_2} \right),\label{e.seconda}
\end{eqnarray}
for some constant $c_n>0$ depending only on $n$.
Let us set for simplicity
$$
\Phi(t)=\max_{\mm_t} \frac{|\nabla A|^2}{g_1g_2}, \qquad \Psi(t) = \min_{\mm_t} |A|^2, \qquad \tilde \sigma=2
\sigma^2\frac{n+2}{3n}
$$
Then we obtain, by the maximum principle,
that $\Phi'(t)<0$ whenever $\Phi(t) > c_n/\tilde\sigma$. We claim that $\Phi(t) \leq c_n/\tilde\sigma$ for all $t<0$. Suppose on the contrary that there is $T_0<0$ such that $\Phi(T_0) > c_n/\tilde\sigma$. Then $\Phi$ is decreasing in $(-\infty,T_0)$. Therefore, we can also find $\dt>0$ such that 
$(1-\dt)\Phi(t) \geq c_n/\tilde \sigma$ on $(-\infty,T_0]$. Then \eqref{e.seconda} implies
$$
\Phi'(t) \leq  -\delta \tilde \sigma \Psi(t) \Phi(t)^2, \qquad t \in (-\infty,T_0),
$$
that is, $(1/\Phi(t))' \geq \delta \tilde\sigma\Psi(t).$
Integrating we obtain
$$
\frac{1}{\Phi(t)} \leq \frac{1}{\Phi(T_0)} - \delta \tilde \sigma \int_{t}^{T_0} \Psi(\tau) d\tau. 
$$
However, this is a contradiction, because the right hand side tends to $-\infty$ as $t \to -\infty$, by assumption \eqref{hyp4b}. This proves that  $\Phi(t)$ remains bounded for all times, which implies the assertion. \qed

\begin{Remark}{\rm An important class of solutions satisfying \eqref{hyp4b} are the ones satisfying the following noncollapsing property of Andrews \cite{A2}: there exists $\alpha>0$ such that for every point $p \in \mm_t$ we can find a ball of radius $\alpha/H(p,t)$ which is enclosed by $\mm_t$ and touches $\mm_t$ from inside at $p$. Such a property is preserved by the flow for any mean convex solution. If we have an ancient convex solutions which is uniformly noncollapsed for $t \in (-\infty,0)$, then at any point we have, using \eqref{eq4},
$$
\frac{1}{|A|^2} \leq \frac{n}{H^2} \leq \frac{n}{\alpha^2} \rho_-(t)^2 \leq \frac{2n^2}{\alpha^2}t,
$$
which implies \eqref{hyp4b}. We recall that for mean convex noncollapsed solutions, not necessarily ancient nor $k$-convex, Haslhofer and Kleiner \cite{HK} have proved a gradient estimate of the form \eqref{grest}, with $C$ also depending on the constant in the noncollapsing condition.}
\end{Remark}

\section{Ancient solutions in the sphere}
\setcounter{equation}{0}

In this paragraph we give some results about ancient solutions immersed in a spherical space form instead of a euclidean space. Therefore, we consider a mean curvature flow $F:\mm \times (-\infty,0) \to \mathbb{S}^{n+1}_R$, where $\mathbb{S}^{n+1}_R$ is the $(n+1)$-dimensional sphere of radius $R>0$ and of sectional curvature $K=1/R^2$. 

We remark that in this case mean curvature flow of a closed hypersurface does not necessarily develop singularities as time increases; for example, minimal submanifolds of the sphere are stationary solutions of the flow. The simplest example is an {\em equator}, i.e. a totally geodesic submanifold of $\mathbb{S}^{n+1}_R$ isometric to $\mathbb{S}^{n}_R$. We call ancient both the solutions which are defined in a maximal interval of the form $(-\infty,T)$, where $T<+\infty$, and the ones defined for all $t \in (-\infty,+\infty)$, which are also called eternal solutions. In the former case, we assume as before $T=0$. There are ancient solutions given by a family of geodesic spheres shrinking around a point, with radius decreasing from $\pi R/2$ to $0$ as $t$ increases. In particular, $\mm_t$ tends to an equator as $t \to -\infty$ and shrinks to a point as $t \to 0$. We call any of these solutions a {\em shrinking spherical cap}. 

\begin{Theorem}
Let $\mm_t$ be an ancient solution of the mean curvature flow in the sphere $\mathbb{S}_R^{n+1}$.
\begin{description}
\item[(a)] If there exists $C>0$ such that $0<|A|^2 \leq C H^2$ for all $t<0$, then $\mm_t$ is a shrinking spherical cap.
\item[(b)] If $\mm_t$ satisfies for all times
\begin{equation}\label{sphpin}
|A|^2-\frac{1}{n-1} H^2 \leq 2K \ \mbox{\rm (if $n \geq 3$)},  \qquad
|A|^2-\frac{3}{4} H^2 \leq \frac {4-\ep}3K \ \mbox{\rm (if $n = 2$)},
\end{equation}
for some $\ep>0$, then $\mm_t$ is either a shrinking spherical cap or an equator.
\end{description}
\end{Theorem}
{\bf Proof} --- (a) Let us again consider the function
$$
f=\frac{|A|^2-\frac{1}{n}H^2}{H^2},
$$
which in this case satisfies, see \cite[Lemma 5.2]{Hu86},
$$
\frac{\partial f}{\partial t}  = 
\Delta f + \frac{2}{H}\langle \nb H, \nb
f \rangle - {\displaystyle \frac{2}{H^4}}|H\nb_i h_{kl}-\nb_i H \,
h_{kl}|^2  - 4nK f. 
$$
If $\mm_t$ is not totally umbilical for all times, there exists $t_1<0$ such that $f \not\equiv 0$ on $\mm_{t_1}$. Applying the maximum principle  we find
$$
0<\max_{\mm_{t_1}} f \leq e^{-4nK(t_1-t)} \max_{\mm_t} f
$$ 
for all $t<t_1$. In particular, we obtain that $\max_{\mm_t} f \to +\infty$ as $t \to -\infty$. On the other hand, the assumption that $|A|^2 \leq C H^2$ implies that $f$ is bounded for all times. The contradiction shows that the solution must be totally umbilical for all times, hence a family of geodesic spheres. We conclude that $\mm_t$ is a shrinking spherical cap.

(b) For fixed $b>0$, we consider the function
$$
\phi_{b}=\frac{|A|^2-\frac{1}{n}H^2}{H^2+b}.
$$
Then, see \cite[Lemma 2.2]{Hu87a}, $\phi_b$ satisfies
\begin{eqnarray}
\frac{\partial \phi_b}{\partial t} & = &
\Delta \phi_b + \displaystyle \frac{4H}{H^2+b}\langle \nb H, \nb
\phi_b \rangle - \frac{2}
{H^2+b}\left(|\nabla A|^2-\left(\phi_b + \frac 1n \right)|\nb H|^2 \right)\nonumber \\ \nonumber \\
& & +2 \phi_b \displaystyle \frac{b|A|^2-2nKH^2-bnK}{H^2+b}.  \label{eqpb}
 \end{eqnarray}
 We consider first the case $n\geq 3$. If we choose
 $$
 b=\frac{3}{2} n(n-1)K
 $$
 then, by \eqref{sphpin}, we find
 \begin{eqnarray} 
 \phi_b & \leq & \frac{\frac{1}{n-1}H^2+2K-\frac{1}{n}H^2}{H^2+b} = \frac{\frac{2}{n(n-1)}H^2+4K}{2H^2+3n(n-1)K}  \nonumber \\
 & < & \frac{4}{3n(n-1)} < \frac{3}{n+2} - \frac 1n. \label{pbunb}
\end{eqnarray}
In view of the well know inequality 
$|\nabla H|^2 \leq \frac{n+2}{3}|\nabla A|^2$, we find that the last term in the first line of \eqref{eqpb} is nonpositive. In addition
 \begin{eqnarray*}
 b|A|^2-2nKH^2-bnK & \leq & \left( \frac{b}{n-1} -2nK \right)H^2 + 2bK-bnK \\
 & = & -\frac {nK}2 H^2 -bK(n-2) < -K (H^2+b).
 \end{eqnarray*}
 We conclude that
 $$
 \frac{\partial \phi_b}{\partial t} \leq 
\Delta \phi_b + \displaystyle \frac{4H}{H^2+b}\langle \nb H, \nb
\phi_b \rangle -2 K \phi_b.
 $$
We now argue as in case (a), finding that if $\phi_b$ is not identically zero then it must become unbounded as $t \to -\infty$, in contradiction with \eqref{pbunb}. The difference with case (a) is that now $H=0$ is allowed, so that a totally umbilical $\mm_t$ can also be an equator.
 
In the case $n=2$,  if we choose $b=\frac{4(4-\ep)}{3}K$, then by \eqref{sphpin} we have
$$
\phi_b+\frac 12 \leq \frac{3H^2 + b-2H^2 }{4(H^2+b)} + \frac 12= \frac 34,
$$
$$
 b|A|^2-4KH^2-2bK  \leq  \left( \frac{3b}{4} -4K \right)H^2 + b\left(\frac{4-\ep}3-2\right)K \leq - \ep (H^2+b).
$$
We can conclude as in the case $n \geq 3$.
\qed

\begin{Remark}{\rm 
Part (a) of the theorem can be applied in particular to convex solutions, where we have $|A|^2<H^2$. More generally, if $\mm_t$ is $k$-convex for some $k=2,\dots,n-1$, we have
$$
0 \leq \ld_1+\dots+\ld_k \leq \ld_1 +(k-1)\ld_k < \ld_1 + n \ld_n,
$$
which implies
$$
|A|^2 \leq n \max_{1 \leq i \leq n} |\ld_i|^2 = n (\max \{ -\ld_1, \ld_n\})^2 \leq n^3 \ld_n^2 < n^3 H^2.
$$
The proof of (a) shows that the result holds under much weaker hypotheses, e.g., if $\mm_t$ satisfies a bound of the form $|A|^2 < e^{B|t|}H^2$ for some $B<4nK$.

The constants in part (b) of the statement are optimal in the case $n \geq 3$. In fact, as observed in \cite{Hu87a}, there exist minimal tori embedded in $\mathbb{S}^{n+1}$ satisfying $|A|^2 < \frac{1}{n-1}H^2+\ep$ with $\ep>0$ arbitrarily small.
}\end{Remark} \medskip

\noindent {\bf Acknowldegments:} 
Gerhard Huisken has been supported by SFB647 of DFG. Carlo Sinestrari has been supported by FIRB-IDEAS project ``Analysis and beyond'', by PRIN 2009 Project ``Viscosity, metric and control theoretic methods in nonlinear PDE's'' and by the group GNAMPA of INdAM (Istituto Nazionale di Alta Matematica).

\end{document}